\documentclass[journal]{IEEEtran}
\pdfoutput=1
\usepackage[latin1,utf8]{inputenx}
\usepackage[dvips]{graphicx}
\usepackage[cmex10]{amsmath}
\usepackage{array}
\usepackage[colorlinks=true,citecolor=blue,linkcolor=blue,urlcolor=blue,plainpages=false]{hyperref}
\usepackage{booktabs, tabularx, multirow, array}
\usepackage{nomencl}
\usepackage{multirow}
\usepackage{comment} 
\usepackage{amsmath}
\usepackage{amsfonts}
\usepackage{amssymb}
\usepackage{graphicx}
\usepackage{amsthm}
\usepackage{eurosym}
\usepackage{algorithm}
\usepackage{tabularx}
\usepackage{cite}
\usepackage{array}
\usepackage{graphicx}
\usepackage{amsfonts,amssymb,amsxtra,balance}
\usepackage{layouts}
\usepackage{accents}
\usepackage{multicol}

\usepackage{enumitem}
 \usepackage[table,xcdraw]{xcolor}

\usepackage{epstopdf}
\usepackage{algpseudocode}
\usepackage{enumitem}

\usepackage{multirow}
\usepackage{dcolumn}
\usepackage{color}

\usepackage{setspace}
\usepackage[export]{adjustbox}
\usepackage{adjustbox}
\usepackage{soul}
\usepackage{booktabs}

\usepackage[cmex10]{amsmath}

\usepackage{tikz}
\usetikzlibrary{shapes,arrows}

\usepackage[textsize=scriptsize,shadow]{todonotes}

\DeclareGraphicsExtensions{.pdf,.png,.jpg}
\graphicspath{{images/}}

\DeclareMathOperator*{\argmin}{arg\,min}

\newcommand\omicron{O}
 
\usepackage{marginnote}
\makenomenclature

\begin{document}
\title{Combining Deep Learning and Optimization for Security-Constrained Optimal Power Flow}

\author{Alexandre~Velloso and Pascal~Van~Hentenryck

\thanks{The authors are affiliated with the School
of Industrial and Systems Engineering, Georgia Institute of Technology, Atlanta, GA 30332, USA. E-mail: avelloso@gatech.edu, pvh@isye.gatech.edu.}
}

\maketitle

\begin{abstract}

The security-constrained optimal power flow (SCOPF) is fundamental in
power systems and connects the automatic primary response (APR) of
synchronized generators with the short-term schedule. Every day, the
SCOPF problem is repeatedly solved for various inputs to determine
robust schedules given a set of contingencies. Unfortunately, the
modeling of APR within the SCOPF problem results in complex
large-scale mixed-integer programs, which are hard to solve. To
address this challenge, leveraging the wealth of available historical
data, this paper proposes a novel approach that combines deep learning
and robust optimization techniques.  Unlike recent machine-learning
applications where the aim is to mitigate the computational burden of
exact solvers, the proposed method predicts directly the SCOPF
implementable solution. Feasibility is enforced in two steps. First,
during training, a Lagrangian dual method penalizes violations of
physical and operations constraints, which are iteratively added as
necessary to the machine-learning model by a
Column-and-Constraint-Generation Algorithm (CCGA). Second, another
different CCGA restores feasibility by finding the closest feasible
solution to the prediction. Experiments on large test cases show that
the method results in significant time reduction for obtaining
feasible solutions with an optimality gap below 0.1\%.

\end{abstract}

\begin{IEEEkeywords}
$\,$ Column-and-constraint generation, decomposition methods, deep learning, neural network, primary response, security-constrained optimal power flow.
\end{IEEEkeywords}


\vspace{-0.35cm}
\section*{Nomenclature}

This section introduces the main notation. Bold symbols are used for
matrices (uppercase) and vectors (lowercase). Additional symbols are
either explained in the context or interpretable by applying
the following general rules: Symbols with superscript
``$j$", ``$k$", or ``$l$" denote new variables, parameters or sets
corresponding to the $j$-th, $k$-th, or $l$-th iteration of the
associated method. Symbols with superscript ``$*$" denote the optimal
value of the associated (iterating) variable. Symbols with superscript
$``t"$ are associated with the data set for the $t$-th past
solve. Dotted symbols are associated with predictors for the
corresponding variable.

\vspace{0.15cm}


\noindent {\bf Sets}

\begin{description} [labelindent=5.8pt ,labelwidth=35pt, labelsep=4pt, leftmargin =44.8pt, style =standard, itemindent=0pt]
	

	\vspace{0.021cm}  \item[$\mathcal{E},\mathcal{E}_s$] Feasibility sets for variables associated with the nominal state and contingent state $s$, respectively.
	
	
	\vspace{0.021cm}  \item[$\mathcal{F}_s$] Feasibility set for primary response variables under contingent state $s$.

	\vspace{0.021cm}  \item[$\mathcal C$, $\mathbb{C}$] Full set and subset of constraints.	
	
	
	\vspace{0.021cm}  \item[$\mathcal G, \mathcal L, \mathcal N$] Sets of generators, transmission lines and buses, respectively. 
	

	\vspace{0.021cm}  \item[$\mathcal S$, $\mathbb{S}$] Full set and subset of contingencies, respectively.
	
	

	\vspace{0.021cm}  \item[$\mathcal{T}$, $\mathbb{T}$] Full set and subset of  past solves, respectively.
	
	
	\vspace{0.021cm}  \item[$\mathbb{U}^+$, $\mathbb{U}^-$] Subsets of line--contingent state pairs.
	
	\vspace{0.021cm}  \item[$\mathcal{Y}_s$] Set of decision variables associated with automatic primary response under contingent state $s$.
		
	

	
	
	
\end{description}

\noindent {\bf Parameters}

\begin{description} [labelindent=5.8pt ,labelwidth=35pt, labelsep=4pt, leftmargin =44.8pt, style =standard, itemindent=0pt]


	\vspace{0.021cm}  \item[$\alpha$, $\rho$] Learning rate and Lagragian dual step size.
	
	\vspace{0.02cm}  \item[$\beta,\beta_1,\beta_c$] Parameters for selecting constraints.
		
	\vspace{0.021cm} \item[$\boldsymbol{\gamma}$] Vector of parameters for primary response.

	\vspace{0.021cm}  \item[$\gamma_i$] Parameter for primary response of generator $i$.

	\vspace{0.021cm}  \item[$\epsilon$] Tolerance for transmission line violation.

	\vspace{0.021cm}  \item[$\boldsymbol{\lambda}, \lambda_c$] Vectors for all Lagrangian multipliers and Lagrangian multipliers for constraint $c$.	
	
	\vspace{0.021cm}  \item[$\boldsymbol{\nu}, \nu_c, \tilde{\nu}_c$] Vector for violations, violation for constraint $c$, and median violation for $c$ among past solves $\mathcal{T}$.

	\vspace{0.021cm}  \item[$\mathbf{A}$, $\mathbf{B}$] Line-bus and Generator-bus incidence matrices.
	
	\vspace{0.021cm}  \item[$\boldsymbol{\omega}$] Vector of weights for deep neural network.

		
	\vspace{0.021cm}  \item[$\mathbf{{d}}$] Vector of nodal net loads.
	
	\vspace{0.021cm}  \item[$\mathbf{e}$] Vector of ones with appropriate dimension.
	
	
	\vspace{0.021cm}  \item[$\mathbf{\overline{f}}$] Vector of line capacities.
	
	\vspace{0.021cm}  \item[$\mathbf{\underline{g}}, \mathbf{\overline{g}}$] Vectors of lower and upper limits for generators.
	
	\vspace{0.021cm}  \item[$\overline{g}_i$] Upper limit for  generator $i$.
		
	\vspace{0.021cm}   \item[$\mathbf{\hat{g}}$] Vector of capacities for generators.

	\vspace{0.021cm}  \item[$\hat{g}_i$] Capacity of generator $i$.
	
	\vspace{0.021cm}  \item[${h(\cdot)}$] Piecewise linear generation costs.
	
	\vspace{0.021cm}  \item[$\mathbf{K}_0$] Matrix of power transfer distribution factors.
	
	\vspace{0.021cm}  \item[$\mathbf{K}_1$] Preprocessed matrix for flow limits.
	
	
	\vspace{0.021cm}  \item[$\mathbf{k}_2, \mathbf{k}_3$] Preprocessed vectors for flow limits.
	
		
	
	
		
	\vspace{0.021cm}  \item[$\mathbf{\overline{r}}$] Vector of primary response limits of generators.
	
	\vspace{0.021cm}  \item[$\overline{r}_i$] Element of $\mathbf{\overline{r}}$ related to generator $i$, given by $\gamma_i\hat{g}_i$.
	
	\vspace{0.021cm}  \item[$\mathbf{S}$] Angle-to-flow matrix.
	

	

	
	
	
	
	
	
	
	
	
\end{description}

\vspace{0.2cm}

\noindent {\bf Nominal-state-related decision variables and vectors}

\begin{description} [labelindent=5.8pt ,labelwidth=35pt, labelsep=4pt, leftmargin =44.8pt, style =standard, itemindent=0pt]

	\vspace{0.02cm} \item[$\mathbf{\boldsymbol{\theta}},\mathbf{f},\mathbf{g}$] Phase angles, line flows, and nominal generation.
	


	\vspace{0.021cm}  \item[$g_i$] Generation of generator $i$ in nominal state.



\end{description}

\vspace{0.08cm}

\noindent {\bf Contingent-state-related decision variables and vectors}

\begin{description}[labelindent=5.8pt ,labelwidth=31pt, labelsep=4pt, leftmargin =40.8pt, style =standard, itemindent=0pt]

	\item[$\mathbf{\boldsymbol{\theta}}_s$] Vector of phase angles under contingent state $s$.



	\vspace{0.021cm}  \item[$\boldsymbol{\tau}_s^+,\boldsymbol{\tau}_s^-$] Vectors of line violation under contingent state $s$.
	
	\vspace{0.021cm}  \item[$\phi$, $s_{\phi}$] Highest line violation and related contingent state.
		\vspace{0.021cm}  \item[$\tilde{\phi}$] Median highest line violation among instances $\mathcal{T}$.
		
	\vspace{0.021cm}  \item[$\mathbf{f}_{s}$] Vector for line flows under contingent state $s$.

	\vspace{0.021cm}  \item[$\mathbf{g}_{s}$] Vector for generation under contingent state $s$.
	
	\vspace{0.021cm}  \item[$\mathbf{g}_{s}^{'}$] Provisional vector for $\mathbf{g}_{s}$.
	
	\vspace{0.021cm}  \item[$g_{s,i}$] Generation of generator $i$ under contingent state $s$.
	
	\vspace{0.021cm}  \item[$g_{s,i}^{'}$] Provisional variable for $g_{s,i}$.	
	
	\vspace{0.021cm}  \item[$n_s$] Global signal under contingent state $s$.


	
	
	\vspace{0.021cm}  \item[$\mathbf{x}_{s}$] Binary vector indicating whether generators reached $\mathbf{\overline{g}}$ under contingency state $s$.
	
	\vspace{0.021cm}  \item[$x_{s,i}$] Element of $\mathbf{x}_{s}$ corresponding to generator $i$.


\end{description}

\section{Introduction}

\subsection{Motivation}

Power systems operations require constant equilibrium between nodal
loads and generation. At the scale of seconds, this balance is
achieved by \emph{Automatic Primary Response} (APR) mechanisms that
govern the \emph{synchronized generators}. For longer time scales,
ranging from a few minutes to hours or even days ahead, this balance
is obtained by solving mathematical optimization problems, as
independent system operators seek consistent and efficient schedules
satisfying complex physical and operational constraints. The need to
solve these optimization problems in a timely manner is driving
intense research about new models and algorithms, both in industry and
academia. In this vein, this work aims at speeding up solution times
of security-constrained optimal power flow (SCOPF) problem \cite{
  alsac1974optimal , bouffard2005umbrella, li2008decomposed,
  capitanescu2011state, wang2016solving,
  dvorkin2016optimizing,velay2019fully} by combining deep learning and
robust optimization methods. The SCOPF is solved by operators every
few minutes for different {sets of bus loads. The high penetration of renewable sources of
energy has increased the frequency of these optimizations}.  The SCOPF
problem considered in this work links the APR to the very short-term
schedule. It is also relevant to mention that the SCOPF problem is
directly or indirectly present in many other power system
applications, including security-constrained unit commitment
\cite{street2010contingency}, transmission switching
\cite{khodaei2012security}, and expansion planning
\cite{moreira2014adjustable}. Thus, a reduction in the computational
burden would allow system operators to introduce important modeling
improvements to many applications.


\subsection{Contextualization and Related Work}

The SCOPF problem determines a least-cost pre-contingency generator
dispatch that allows for feasible points of operation for a set of
\emph{contingencies}, e.g., individual failures of main lines and/or
generators. The SCOPF problem may refer to the \emph{corrective} case
\cite{wang2016solving} where re-scheduling is deemed possible and to
the \emph{preventive} case where no re-dispatch occurs
\cite{li2008decomposed, dvorkin2016optimizing, velloso2019exact},
i.e., the system must be able to achieve a feasible steady-state point
without a new schedule. A valuable review of the SCOPF problem and
solution methods is available in
\cite{capitanescu2011state}. Interesting discussions about credible
contingencies, reserve requirements, security criteria, and regulation
for reserves can be found in
\cite{bouffard2005umbrella,dvorkin2016optimizing,zhou2016survey} and
the references therein. Without loss of generality, the $N-1$ security
criterion for generators is adopted in this paper, i.e., the system
must operate under the loss of any single generator.

The SCOPF is a nonlinear and nonconvex problem based on the AC {optimal power flow (OPF)}
equations. Extensive reviews can be found in \cite{acscopf} and
\cite{capit2016critical}. The DC formulation of the SCOPF has been
widely used in both academia and industry \cite{bouffard2005umbrella,
  li2008decomposed, dvorkin2016optimizing,velay2019fully}. Interesting
discussions regarding the quality of approximations and relaxations of
the OPF problem can be found in \cite{li2007dcopf,coffrin2014linear,
  coffrin2015qc,eldridge2017marginal}. The DC-SCOPF can also be used
to improve AC-SCOPF approaches \cite{2012dcexploiting}. It is not
within the scope of this work to discuss the quality of aforementioned
approximations or relaxations to the optimal power flow. Instead, this
research offers a new approach that improves current industry
practices, which is still strongly based on DC-SCOPF.

The APR of synchronized generators is essential for stability. These
generators respond automatically to frequency variations, caused by
power imbalances {for instance}, by adjusting their power outputs until
frequency is normalized and the power balance is
restored. Unfortunately, the APR deployment, which is bounded by
generators limits only, may result in transmission line overloads
\cite{dvorkin2016optimizing, lannoye2014transmission,
  velloso2019exact}. Therefore, this work co-optimizes the APR of
synchronized generators within the (preventive) SCOPF problem. Even
though the APR behavior is nonlinear, linear approximations are used
in practice \cite{kundur1994power}. In
\cite{restrepo2005unit,dvorkin2016optimizing,
  velay2019fully,karoui2010modeling, velloso2019exact}, the APR is
modeled by a single variable representing frequency drop (or power
loss) for each contingency state and by a participation factor for
each generator.

The DC-SCOPF problem with APR is referred to as the SCOPF problem for
conciseness in this paper. It admits an exact \emph{extensive
  formulation} (with all variables and constraints for nominal state
and contingency states) as a mixed-integer linear program (MILP)
\cite{dvorkin2016optimizing, velloso2019exact}. Nevertheless, this
formulation is generally very large because the APR constraints
require binary variables for each generator and for each contingency
state to determine whether generators are producing according to the
linear response model or are at their limits
\cite{dvorkin2016optimizing, velloso2019exact}. Thus, the number of
binary variables increases quadratically with the number of
generators, which makes the extensive form of the SCOPF
impractical. Better modeling strategies and decomposition schemes are
required.

The robust optimization framework has been widely applied in power
systems due to its interesting tradeoff between modeling capability
and tractability. See \cite{ShabbirBook} for a review of robust
optimization applications in power systems. The SCOPF problem can be
modeled as a two-stage robust optimization or adaptive robust
optimization (ARO) and tackled by two main decomposition methods:
Benders decomposition \cite{bertsimas2013adaptive} and the
{column-and-constraint-generation algorithm (CCGA)}
\cite{BoZeng2011}. Both approaches rely on iterative procedures that
solve a \emph{master problem} and \emph{subproblems}.  The master
problem is basically the nominal OPF problem with additional
cuts/constraints and variables representing feasibility or optimality
information on the subproblems. Whereas Benders decomposition provides
dual information about the subproblems through valid cuts restricting
the master problem, the CCGA adds primal contraints and variables from
the subproblems to the master problem.

Unfortunately, the SCOPF problem is not suitable for a traditional
Benders decomposition since the subproblems are nonconvex due to the
APR constraints (which feature binary variables). Despite such
challenges, inspired by \cite{bertsimas2013adaptive}, an interesting
heuristic method was proposed in \cite{dvorkin2016optimizing} but it
does not guarantee optimality. In contrast, the CCGA algorithm
proposed in \cite{velloso2019exact} is an exact solution method which
was used to produce optimal solutions to power network with more than
2,000 buses. Notwithstanding aforementioned contributions, both
approaches still require significant computational effort to obtain
near-optimal solutions.

Machine learning (ML) approaches have been advocated to address the
computational burden associated with the hard and repetitive
optimization problems in the power sector, given the large amounts of
historical data (i.e., past solutions). Initial attempts date back to
the early 1990s, when, for example, artificial neural networks were
applied to predict the on/off decisions of generators for an unit
commitment problem \cite{sasaki1992solution,ouyang1992hybrid}.
More recently, {ML} was used to identify partial
warm-start solutions and/or constraints that can be omitted, and to
determine affine subspaces where the optimal solution is likely to
lie \cite{xavier2019learning}.  Artificial neural network and
decision tree regression were also used to learn sets of
high-priority lines to consider for transmission switching
\cite{yang2019line}, while the k-nearest neighbors approach was used
to select previously optimized topologies directly from data
\cite{johnson2020k}. As for the security and reliability aspects of
the network, the security-boundary detection was modeled with a
neural network to simplify stability constraints for the optimal
power flow \cite{gutierrez2010neural}, while decision trees were
{applied to determine security boundaries (regions) for controllable
variables for a coupled natural gas and electricity system
\cite{costa2016decision}}. Machine learning was also applied
for identifying the relevant sets of active constraints for the OPF
problem \cite{misra2018learning}.

Unlike these applications, where the main purpose of machine learning
is to enhance the solver performance by classifying sets, eliminating
constraints, and/or by modeling specific parts of the problems, the
machine-learning approach in \cite{fioretto2019predicting} directly
predicts the generator dispatch for the OPF by combining deep learning
and Lagrangian duality. This approach produces significant
computational gains but is not directly applicable to the SCOPF
problem which features an impractical number of variables and
constraints. This work remedies this limitation.

\subsection{Contributions}

The paper assumes the existence of historical {SCOPF} data, i.e., pairs
of inputs and outputs
\cite{misra2018learning,xavier2019learning,yang2019line,johnson2020k,fioretto2019predicting}
The proposed approach uses a deep neural network (DNN) to approximate
the mapping between loads and optimal generator dispatches. To capture
the physical, operational, and APR constraints, the paper applies the
Lagrangian dual scheme of \cite{fioretto2019predicting} that penalizes
constraint violations at training time. Moreover, to ensure
computational tractability, the training process, {labeled as CCGA-DNN,} mimics a dedicated
CCGA algorithm that iteratively adds new constraints for a few
critical contingencies. In these constraints, an approximation for the
post-contingency generation is adopted to keep the size of the DNN
small. The resulting DNN provides high-quality approximations to the
SCOPF in milliseconds and can be used to seed another dedicated CCGA
to find the nearest feasible solution to the prediction. The resulting
approach may bring two orders of magnitude improvement in efficiency
compared to the original CCGA algorithm. 

In summary, the contributions can be summarized as follows: i) a novel
DNN that maps a load profile onto a high-quality approximation of the
SCOPF problem, ii) a new training procedure{, the CCGA-DNN,} that mimics a CCGA, where
the master optimization problem is replaced by a DNN {prediction}, iii)
an approximation for the post-contingency generation which keeps the
DNN size small, and iv) an dedicated CCGA algorithm seeded with the
DNN evaluation to obtain high-quality feasible solutions fast. Of
particular interest is the tight combination of machine learning and
optimization proposed by the approach. 
  
\subsection{Paper Organization}

This paper is organized as follows. Section \ref{sec.SCOPF} introduces
the SCOPF problem. Section \ref{sec.Meth} presents the properties of
the SCOPF problem and the CCGA for SCOPF. Section \ref{sec.DNN}
introduces the {deep learning }models in stepwise refinements. Section
\ref{sec.restaur} describes the CCGA for feasibility recovery.
Section \ref{sec.results} reports the case studies and the numerical
experiments and Section \ref{sec.conclusion} concludes the paper.

\section{The SCOPF Problem}
\label{sec.SCOPF}


\subsection{The Power Flow Constraints}

The SCOPF formulation uses traditional security-constrained DC power
flow constraints over the vectors for generation $\mathbf{g}$, flows
$\mathbf{f}$, and phase angles $\mathbf{\boldsymbol{\theta}}$. In matrix
notations, these constraints are represented as follows:
\vspace{-0.5cm}
\begin{multicols}{2}
	\begin{flalign}
	&\mathbf{A} \mathbf{f} + \mathbf{B} \mathbf{g} = \mathbf{{d}} \label{eq.Master.EnergyBalance}& \hspace{-1cm} \\
	&\mathbf{f} = \mathbf{S}  \mathbf{\boldsymbol{\theta}}  \label{eq.sec.Kirc}& \hspace{-1cm} \\ 
	&-\mathbf{\overline{f}} \leq \mathbf{f} \leq \mathbf{\overline{f}} \label{eq.Master.PFLimit}& \hspace{-1cm} \\
	&\mathbf{\underline{g}}\leq\mathbf{g} \leq \mathbf{\overline{g}}\label{eq.Master.GenCap}& \hspace{-1cm}
	\end{flalign}
	\columnbreak\\
\begin{flalign}
	&\hspace{-1cm} \mathbf{A} \mathbf{f}_s + \mathbf{B} \mathbf{g}_s = \mathbf{{d}} & \forall s \in \mathcal S \label{eq.Conting.EnergyBalance}\\
	&\hspace{-1cm} \mathbf{f}_s = \mathbf{S}  \mathbf{\boldsymbol{\theta}}_s & \forall s \in \mathcal S \label{eq.sec.Kirc.scenario}\\
	&\hspace{-1cm} -\mathbf{\overline{f}} \leq \mathbf{f}_s \leq \mathbf{\overline{f}} \label{eq.PFLimit.scenario} & \forall s \in \mathcal S\\
	&\hspace{-1cm} \mathbf{g}_s \leq \mathbf{\overline{g}} &\forall s \in \mathcal S\label{eq.Scenario.GenCap}
\end{flalign}
\end{multicols}
\vspace{-0.50cm}

\noindent
Equations \eqref{eq.Master.EnergyBalance}--\eqref{eq.Master.GenCap}
model the DC power flow in pre-contingency state and capture the nodal
power balance \eqref{eq.Master.EnergyBalance}, Kirchhoff's second law
\eqref{eq.sec.Kirc}, transmission line limits
\eqref{eq.Master.PFLimit}, and generator limits
\eqref{eq.Master.GenCap}. Analogously, equations
\eqref{eq.Conting.EnergyBalance}--\eqref{eq.Scenario.GenCap} model the
power flow for each post-contingency state $s$. The bounds
$\mathbf{\overline{g}}$ in \eqref{eq.Master.GenCap} and
\eqref{eq.Scenario.GenCap} may be different from capacity
$\mathbf{\hat{g}}$ due to commitment and/or operational constraints.


\subsection{Automatic Primary Response}

The APR is modeled as in
\cite{velloso2019exact,dvorkin2016optimizing,velay2019fully}: under
contingent state $s$, a global variable $n_s$ is used to mimic the
level of system response required for adjusting the power
imbalance. The APR of generator $i$ under contingency $s$,
$g_{s,i}-g_i$, is proportional to its capacity $\hat{g}_i$ and to the
parameter $\gamma_i$ associated with the droop coefficient, i.e.,
\begin{align}
&g_{s,i}= \min \{g_i + n_s \gamma_i \,\hat{g}_i,\;\overline{g}_i\} & &\forall i \in \mathcal{G}, \forall s \in \mathcal{S}, i\ne s
\label{disjunction_1}\\
&g_{s,s}=0 & & \forall s \in \mathcal{S} \label{n-1constraint}.
\end{align}

\noindent
These equations are nonconvex and can be linearized by introducing binary variables $x_{s,i}$ to denote
whether generator $i$ in scenario $s$ is not at its limit, i.e.,
\begin{flalign}
&  |g_{s,i}-g_i- n_s\gamma_i \hat{g}_i| \leq\label{disj01} \overline{g}_i(1-x_{s,i}) \hspace{-0.5cm}  & \forall i \in \mathcal G, s \in \mathcal S, i \ne s \\
& g_i + n_s \gamma_i \,\hat{g}_i \geq \overline{g}_i (1-x_{s,i}) & \forall i \in \mathcal G, s \in \mathcal S, i \ne s \\
& g_{s,i} \geq \overline{g}_i (1-x_{s,i}) & \forall i \in \mathcal G, s \in \mathcal S, i \ne s \\
& n_s \in [0,1] & \forall s \in \mathcal S \\
& x_{s,i} \in \{0,1\}& \forall i \in \mathcal G, s \in \mathcal S \label{disj05}\\
& g_{s,s}=0 &  \forall s \in \mathcal{S} \label{disj06}.
\end{flalign}
\subsection{Extensive Formulation for the SCOPF Problem}
\label{sec.EF}

\noindent
The extensive formulation for the SCOPF problem using variables for generation, flows, and phase angles is as follows:
\begin{align}
\min_{ \mathbf{\boldsymbol{\theta}},\mathbf{f},\mathbf{g},[ \mathbf{\boldsymbol{\theta}}_s,\mathbf{f}_{s},\mathbf{g}_{s},n_s,\mathbf{x}_s]_{s\in\mathcal{S}}}
\quad&  h(\mathbf{g})     \label{obj.long.func}\\
\text{s.t.:}\quad &          \eqref{eq.Master.EnergyBalance}-\eqref{eq.Master.GenCap} \label{const.long.set01} \\
&            \eqref{eq.Conting.EnergyBalance}-\eqref{disj06} \label{const.long.set02} &\forall s \in \mathcal S.
\end{align}

\noindent
Using power transfer distribution factors (PTDF), constraints
\eqref{eq.Master.EnergyBalance}--\eqref{eq.Scenario.GenCap} can be
replaced by the following constraints:
\vspace{-0.65cm}
\begin{multicols}{2}
	\begin{flalign}
    &\mathbf{e}^{\top} \mathbf{g}   \label{PTDF_nominal_first} = \mathbf{e}^{\top}\mathbf{d} & \hspace{-0.76cm}\\
    &|\mathbf{K}_0 (\mathbf{{d}} - \mathbf{B} \mathbf{g} )|  \leq \mathbf{\overline{f}} & \hspace{-0.76cm}\label{PTDF_bound}\\
    &\eqref{eq.Master.GenCap}& \hspace{-0.76cm}\label{eq.PTDFlast_nominal}
	\end{flalign}
	\columnbreak\\
\begin{flalign}
    &\hspace{-0.80cm}\mathbf{e}^{\top} \mathbf{g}_s   \label{valid.pstcont} = \mathbf{e}^{\top}\mathbf{d} & &\forall s \in \mathcal S\\
    &\hspace{-0.80cm}|\mathbf{K}_0 (\mathbf{{d}} - \mathbf{B} \mathbf{g}_s )| \label{preCut_01} \leq \mathbf{\overline{f}} & &\forall s \in \mathcal S \\
    &\hspace{-0.80cm}\eqref{eq.Scenario.GenCap}& &\label{eq.lastPTDFlast}\forall s \in \mathcal S
\end{flalign}
\end{multicols}
\vspace{-0.65cm}

\noindent
Constraints \eqref{PTDF_nominal_first}--\eqref{eq.lastPTDFlast} that
involve the PTDF matrix $\mathbf{K}_0$ are from
\cite{ardakani2013identification}. The total demand balance for the
nominal and contingent states are enforced by
\eqref{PTDF_nominal_first} and \eqref{valid.pstcont} respectively. In
constraints \eqref{PTDF_bound} and \eqref{preCut_01}, the PTDF matrix
translates the power injected by each generator at its bus into its
contribution to the flow of each line. These constraints also bound
the flows from above and below. Observe that $\mathbf{g}$ and
$\mathbf{g}_s$ are the only variables in this formulation.

For conciseness, denote the power flow constraints
\eqref{PTDF_nominal_first}--\eqref{eq.PTDFlast_nominal} and
\eqref{valid.pstcont}--\eqref{eq.lastPTDFlast} by
$\mathbf{g}\in\mathcal{E}$ and $\mathbf{g}_s\in\mathcal{E}_s$
respectively. Similarly, denote the APR constraints
\eqref{disj01}--\eqref{disj06} by
$\mathcal{Y}_s = [\mathbf{g},\mathbf{g}_s, \mathbf{x}_s, n_s ] \in \mathcal{F}_s$.
The extensive SCOPF formulation then becomes
\begin{flalign}
\min_{ \mathbf{g},[ \mathbf{g}_{s},\mathbf{x}_s, n_s]_{s\in\mathcal{S}}}
\quad&  {h(\mathbf{g})}     \label{obj.func}\\
\text{s.t.:}\quad &          \mathbf{g} \in \mathcal{E} \label{const.set01} \\
&            \mathbf{g}_s \in \mathcal{E}_s \label{const.set02} &\forall s \in \mathcal S \\
&            \mathcal{Y}_s \in \mathcal{F}_s \label{const.set03} &\forall s \in \mathcal S
\end{flalign}

\noindent
Note that the number of binary variables above grows quadratically
with the number of generators. Hence, solving
\eqref{obj.func}--\eqref{const.set03} becomes impractical for
large-scale systems.

\section{SCOPF Properties and CCGA}
\label{sec.Meth}

This section introduces key properties of the SCOPF problem and
summarizes the CCGA proposed in \cite{velloso2019exact}. These
properties are necessary for the CCGA and the ML
models. The CCGA serves both as a benchmark for evaluation and is used
as part of the feasibility recovery scheme proposed in Section
\ref{sec.restaur}.
\vspace{0.1cm}

\noindent\textbf{Property 1:} \emph{For $s\in\mathcal{S}$, given
  values $\mathbf{g}^{*}$ and $n_s^{*}$ for $\mathbf{g}$ and $n_s$,
  there exists a unique value $\mathbf{g}^{*}_s$ for $\mathbf{g}_s$
that can be computed directly using constraints
\eqref{disj01}--\eqref{disj06}}.

\vspace{0.1cm}
\noindent \textbf{Property 2:} \emph{Consider $s\in\mathcal{S}$ and a
  value $\mathbf{g}^{*}$ for $\mathbf{g}$. If there exists a value
  $n_s^{*}$ for $n_s$ that admits a feasible solution to constraints
  \eqref{disj01}--\eqref{disj06} and \eqref{valid.pstcont}, then this value
  $n_s^{*}$ is unique and can be computed by a simple bisection method
  \cite{velloso2019exact}.}
\vspace{0.1cm}

\noindent Property 2 holds since, for a given $\mathbf{g}^{*}$, each
component of $\mathbf{g}_s$ is continuous and monotone with respect to
$n_s$. Hence the value $n_s^*$ and its associated vector
$\mathbf{g}_s^{*}$ that satisfy constraint \eqref{valid.pstcont} can
be found by a simple bisection search over $n_s$.
\vspace{0.1cm}

\noindent \textbf{Property 3:} \emph{Constraint \eqref{preCut_01} can be formulated as
\begin{align}
&  \mathbf{K}_1 \mathbf{g}_s + \mathbf{k}_2\geq \mathbf{0} \label{Cut_01}\\
& \mathbf{K}_1 \mathbf{g}_s + \mathbf{k}_3\geq \mathbf{0}  \label{Cut_02},
\end{align}
\noindent using matrix operations to obtain $\mathbf{K}_1$, $\mathbf{k}_2$, and $\mathbf{k}_3$.}
\vspace{0.1cm}

\noindent Note that each row of $\mathbf{K}_1$ and each element of
$\mathbf{k}_2$ and $\mathbf{k}_3$ are associated with a specific
transmission line. Therefore, for each $s\in\mathcal{S}$ and for each
line, the (positive and negative) violation of the thermal limit of
the line can be obtained by inspecting \eqref{Cut_01}--\eqref{Cut_02}
for the proposed value $\mathbf{g}_s^{(*)}$.

\subsection{The Column and Constraint Generation Algorithm}

The CCGA, which relies on the above properties, alternates between
solving a master problem to obtain a nominal schedule $\mathbf{g}$ and
a bisection method to obtain the state variables of each contingency.
The master problem is specified as follows:

\begin{flalign}
 \min_{\mathbf{g},[\mathbf{g}_{s}^{'}]_{s\in\mathcal{S}},[\mathbf{x}_{s},n_s]_{s\in\mathbb{S}}}
\quad&   {h(\mathbf{g})}   \label{master.ccga} \\
\text{s.t.:}\quad &     \mathbf{g} \in \mathcal{E} \label{const.nom} \\
& \mathbf{g}_{s}^{'} -\mathbf{g} \leq  \mathbf{\overline{r}} & \forall s \in \mathcal S \label{RampingCCGA}\\
& \eqref{eq.Scenario.GenCap}, \eqref{valid.pstcont}, \eqref{disj06}& \forall s \in \mathcal{S}\label{semifinal.ccga}\\
& \mathcal{Y}_s \in \mathcal{F}_s& \forall s \in \mathbb{S}\label{final.ccga}\\
&  \mathbf{K}^{l}_{1} \mathbf{g}_{s}^{'} + \mathbf{k}^{l}_{2}\geq \mathbf{0} \label{Cut_01.ccga}& \forall (l,s) \in \mathbb{U}^+\\
& \mathbf{K}^{l}_1 \mathbf{g}_{s}^{'} + \mathbf{k}^{l}_3\geq \mathbf{0} \label{Cut_02.ccga}& \forall (l,s) \in \mathbb{U}^-
\end{flalign}

\noindent
Constraints \eqref{master.ccga}--\eqref{Cut_02.ccga} uses variables
$\mathbf{g}_{s}^{'}$ to denote a ``guess'' for the post-contingency
generation: the actual vector $\mathbf{g}_{s}$ is not determined by
the master problem but by the aforementioned bisection
method. Constraint \eqref{const.nom} enforces the nominal state
constraints.  Constraint \eqref{RampingCCGA} imposes a valid bound for
post-contingency generation. For all contingencies, constraint
\eqref{semifinal.ccga} enforces the generation capacity
\eqref{eq.Scenario.GenCap}, total demand satisfaction
\eqref{valid.pstcont}, and the absence of generation for a failed
generator \eqref{disj06}. The APR is enforced ``on-demand'' in
\eqref{final.ccga} for a reduced set of contingent states
$\mathbb{S}$. Initially, $\mathbb{S}=\emptyset$. Inequalities
\eqref{Cut_01.ccga}--\eqref{Cut_02.ccga} are also the ``on-demand''
versions of \eqref{Cut_01}--\eqref{Cut_02} for (a few) pairs of
transmission lines and contingencies. Initially, 
$\mathbb{U}^+$ and $\mathbb{U}^-$ are empty sets. 


The CCGA algorithm is specified in Algorithm \ref{Alg.CCGA_feas}.  At
iteration $j$, the master problem
\eqref{master.ccga}--\eqref{Cut_02.ccga} computes $\mathbf{g}^j$.  The
bisection method then determines the contingent state variables $[
  \mathbf{g}_{s}^j,\mathbf{x}_s^j, n_s^j]_{s\in\mathcal{S}}$. The
vectors $\boldsymbol{\tau}_s^+$ and $\boldsymbol{\tau}_s^-$ of
positive numbers represent the positive and negative violations of
transmission lines for contingent state $s$: they are calculated for
all $s\in\mathcal{S}$ by inspecting constraints
\eqref{Cut_01}--\eqref{Cut_02} for
$[\mathbf{g}_{s}^j]_{s\in\mathcal{S}}$. The algorithm then computes
the highest single line violation $\phi$ among all contingent states
and uses $s_{\phi}$ to denote the contingent state associated with
$\phi$. The pairs lines/contingencies featuring violations above a
predefined threshold $\beta$ are added to the master problem by
updating sets $\mathbb{U}^+$ and $\mathbb{U}^-$. Likewise, $s_{\phi}$
is added to $\mathbb{S}$. As a result, the variables and APR
constraints associated with $s_{\phi}$ are added to the master problem
during the next iteration. The CCGA terminates when $\phi<\epsilon$,
where $\epsilon$ is the tolerance for line violation.

\begin{algorithm}[t]
	{\footnotesize \caption{{\small CCGA\label{Alg.CCGA_feas} }}
		\begin{algorithmic}[1]
			\State{Initialization: $j\leftarrow0, \mathbb{S}\leftarrow\emptyset$, $\mathbb{U}^+\leftarrow\emptyset$, $\,\mathbb{U}^-\leftarrow\emptyset$}
			\For{$j=0, 1, \ldots$}
			\State{solve \eqref{master.ccga}--\eqref{final.ccga} to obtain $\mathbf{g}^{j}$}
			\State{$n^j_s \leftarrow$ apply the bisection method on all $s\in\mathcal{S}$}
			\State{$\mathbf{g}^{j}_s\leftarrow$ enforce \eqref{disj01}--\eqref{disj06} on all $s\in\mathcal{S}$}
			\State{$\boldsymbol{\tau}_s^-$, $\boldsymbol{\tau}_s^+\leftarrow$ get the line violations of $\mathbf{g}_s^{j}$ using \eqref{Cut_01}--\eqref{Cut_02} for all $s\in\mathcal{S}$ }
			\State{$\phi\leftarrow$ compute the highest line violation among all $s\in\mathcal{S}$}
			\State{$s_{\phi}\leftarrow$ select the contingent state associated with $\phi$}
			\State{$\mathbb{S}\leftarrow\mathbb{S}\cup \{s_{\phi}\}$}
			\State{$\mathbb{U}^+\leftarrow\mathbb{U}^+\cup \{ \,(l, s) \, | \, \boldsymbol{\tau}^+_{s}[l] > \beta\}$}
			\State{$\mathbb{U}^-\leftarrow\mathbb{U}^-\cup \{ \,(l, s) \, | \, \boldsymbol{\tau}^-_{s}[l] > \beta\}$}
			\State{\textbf{BREAK} if $\phi\leq\epsilon$.}
			\EndFor
		\end{algorithmic}
	} 
\end{algorithm}

The following result from \cite{velloso2019exact} ensures the
correctness of CCGA: It shows that a solution to the master problem
produces a nominal generation for which there exists a solution to
each contingency that satisfies the APR and total demand
constraints.  Since the CCGA adds at least one violated line
constraint and, possibly, a set of violated APR constraints for one
contingency to the master problem at each iteration, it is guaranteed
to converge after a finite number of iterations.

\vspace{0.1cm}
\noindent \textbf{Theorem 1:} \emph{For each solution $\mathbf{g}^{*}$
  to the master problem, there exist values $n^{*}_s$ and
  $\mathbf{g}^{*}_s$ that satisfy the demand constraint
  $\mathbf{e}^{\top}\mathbf{g}^{*}_s = \mathbf{e}^{\top}\mathbf{d}$
  and the APR constraints \eqref{disj01}--\eqref{disj06}} for each
contingency $s$.

\noindent \textbf{Proof}: By \eqref{eq.Scenario.GenCap} and
\eqref{RampingCCGA}, $g_{s,i}^{'}\leq \min{\{\overline{g}_i, g_i +
  \gamma_i \hat{g}_i\}}$ for each $i$ and $s$, where $g_{s,i}^{'}$ is
the $i$-th element of $\mathbf{g}_{s}^{'}$. When $n_s=0$,
$\mathbf{g}_{s}=\mathbf{g}$, except for $g_{s,s}=0$. When $n_s=1$,
$g_{s,i} = \min{\{\overline{g}_i, g_i + \gamma_i \hat{g}_i\}} \geq
g_{s,i}^{'}$ for each $i$ and $s$, with $i\ne s$. Since, by
\eqref{valid.pstcont}, $\mathbf{g}_{s}^{'}$ meets the global demand,
$\mathbf{e}^{\top} \mathbf{g}_s \geq \mathbf{e}^{\top}
\mathbf{g}_s^{'} = \mathbf{e}^{\top}\mathbf{d}$ when $n_s=1$. By the
monotonicity and continuity of $g_{s,i}$ with respect to $n_s$ (for a
given $g_i$), there is a value $n_s^{*}$ whose associated
$\mathbf{g}_{s}^{*}$ satisfies the demand constraint in
\eqref{valid.pstcont} and preserves the APR constraints $\square$.

\section{Deep Neural Networks for SCOPF}
\label{sec.DNN}

This section describes the use of supervised learning to obtain DNNs that map a load vector into a solution of the SCOPF
problem. A DNN consists of many layers, where the
input for each layer is typically the output of the previous layer
\cite{lecun2015deep}. This work uses fully-connected DNNs.

\subsection{Specification of the Learning Problem}
\label{sec.specLearn}

For didactic purposes, the specification of the learning problem
uses the extensive formulation \eqref{obj.func}--\eqref{const.set03}.
The training data is a collection of instances of the form
\[
\{\mathbf{d}^t ; \mathbf{g}^t, [\mathbf{g}^t_s, n^t_s, \mathbf{x}^t_s]_{s\in\mathcal{S}}\}_{t\in\mathcal{T}}
\]
where $(\mathbf{g}^t, [\mathbf{g}^t_s, n^t_s,
  \mathbf{x}^t_s]_{s\in\mathcal{S}})$ is the optimal solution (ground
truth) to the SCOPF problem for input $\mathbf{d}^t$. The DNN is a
parametric function \omicron$[\boldsymbol{\omega}](\cdot)$ whose
parameters are the network $\boldsymbol{\omega}$: It maps a load
vector $\mathbf{d}$ into an approximation
\omicron$[\boldsymbol{\omega}](\mathbf{d}) =
\{\mathbf{\dot{g}},[\dot{n}_s,\mathbf{\dot{x}}_s,\mathbf{\dot{g}}_s]_{s\in\mathcal{S}}\}$
of the optimal solution to the SCOPF problem for load
$\mathbf{d}$. The goal of the machine-learning training to find the
optimal weights $\boldsymbol{\omega}^*$, i.e.,
\begin{flalign}
\boldsymbol{\omega}^* = \argmin_{ \boldsymbol{\omega}}
& \sum_{t\in\mathcal{T}} \mathbb{L}^t_0(\mathbf{\dot{g}}^t) + \sum_{s\in\mathcal{S}} \mathbb{L}^t_s(\mathbf{\dot{g}}^t_s,\mathbf{\dot{x}}^t_s,\dot{n}^t_s)
& \hspace{-2cm} \label{ML.obj.func}\\
\text{s.t.: } & \omicron[\boldsymbol{\omega}](\mathbf{d}^t) = (\mathbf{\dot{g}}^t,[\mathbf{\dot{g}}^t_s, \dot{n}^t_s,\mathbf{\dot{x}}^t_s]) & \hspace{-2cm}  \forall t \in \mathcal T  \label{ML.def}\\
& \mathbf{\dot{g}}^t \in \mathcal{E}^t& \hspace{-2cm}  \forall t \in \mathcal T  \label{ML.const.set01}\\
&\mathbf{\dot{g}}^t_s \in \mathcal{E}^t_s& \hspace{-2cm} \forall t \in \mathcal T, \forall s \in \mathcal S & \label{ML.const.set02}\\
&\dot{\mathcal{Y}}^t_s \in \mathcal{F}^t_s& \hspace{-2cm} \forall t \in \mathcal T, \forall s \in \mathcal S & \label{ML.const.set03}
\end{flalign}
\noindent
where the loss functions are defined as 
\begin{flalign*}
& \mathbb{L}^t_0(\mathbf{\dot{g}}^t) = ||\mathbf{g}^t-\mathbf{\dot{g}}^t||_{2} \\
& \mathbb{L}^t_s(\mathbf{\dot{g}}^t_s,\mathbf{\dot{x}}^t_s,\dot{n}^t_s) = ||\mathbf{g}^t_s-\mathbf{\dot{g}}^t_s||_{2} + ||\mathbf{x}^t_s-\mathbf{\dot{x}}^t_s||_{2}+ ||n^t_s-\dot{n}^t_s||_{2} 
\end{flalign*}
and minimize the distance between the prediction and the ground
truth. There are two difficulties in this learning problem: the large
number of scenarios, variables, and constraints, and the satisfaction
of constraints \eqref{ML.const.set01}--\eqref{ML.const.set03}. This
section examines possible approaches.

\subsection{The Baseline Model}

The baseline model is a parsimonious approach that disregards
constraints \eqref{ML.const.set01}--\eqref{ML.const.set03} and
predicts the nominal generation only, i.e.,
\begin{flalign*}
& \boldsymbol{\omega}^* = & \hspace{-2.0cm} \argmin_{ \boldsymbol{\omega}} & \sum_{t\in\mathcal{T}} \mathbb{L}^t_0(\mathbf{\dot{g}}^t)  \\
& & \hspace{-2.0cm} \text{s.t.: }& \omicron[\boldsymbol{\omega}](\mathbf{d}^t) = (\mathbf{\dot{g}}^t) & \forall t \in \mathcal T 
\end{flalign*}
It uses a DDN model with 5 linear layers, interspersed with 5
nonlinear layers that use the softplus activation function. The sizes
of the input and output of each layer are linearly parameterized by
$|\mathcal{L}|$ and $|\mathcal{G}|$. A high-level algebraic description of layers of the DNN follows:
\begin{flalign*}
\qquad&\mathbf{l}_i = \gamma ( \mathbf{W}_i \mathbf{l}_{i-1} + \mathbf{b}_i),&&\quad \text{for each layer} \;\mathbf{l}_i\qquad\\
&\mathbf{l}_1 = \gamma ( \mathbf{W}_1 \mathbf{d} + \mathbf{b}_1)
\end{flalign*}
The elements of vector $\boldsymbol{\omega}$ are rearranged as
matrices $\mathbf{W}_i$ and vector of biases $\mathbf{b}_i$. Note that
the demand vector $\mathbf{d}$ is the input for the first layer. The
symbol $\gamma$ denotes a nonlinear activation function.
Unfortunately, training the baseline model tends to produce predictors
violating the problem constraints
\cite{misra2018learning,fioretto2019predicting}.

\subsection{A Lagrangian Dual Model for Nominal Constraints}
\label{sec.lagdual}

This section extends the baseline model to include constraints on the
nominal state \eqref{ML.const.set01}. Constraints
\eqref{ML.const.set02}--\eqref{ML.const.set03} on the contingency
cases are not considered in the model. To capture physical and
operational constraints, the training of the DNN adopts the Lagrangian
dual approach from \cite{fioretto2019predicting}.

The Lagrangian dual approach relies on the concept of constraint
violations. The violations of a constraint $f(x) = 0$ is given by
$|f(x)|$, while the violations of $f(x) \geq 0$ are specified by
$\max(0,-f(x))$. Although these expressions are not differentiable,
they admit subgradients. Let $\mathbb{C}$ represent the set of nominal
constraints and $\nu_{c}(\mathbf{g})$ be the violations of constraint
$c$ for generation dispatch $\mathbf{g}$. The Lagrangian dual approach
introduces a term $\lambda_c \nu^t_{c}(\mathbf{g}^t)$ in the objective
function for each $c \in \mathbb{V}$ and each $t \in \mathcal{T}$,
where $\lambda_c$ is a Lagrangian multiplier. The
optimization problem then becomes
\begin{flalign}
& LR(\boldsymbol{\lambda}) = & \hspace{-0.5cm} \min_{ \boldsymbol{\omega}} & \quad \sum_{t\in\mathcal{T}} (\mathbb{L}^t_0(\mathbf{\dot{g}}^t)  + \sum_{c\in\mathcal{C}} \lambda_c \nu^t_{c}(\mathbf{\dot{g}}^t)) \label{abst.min_LR} \\
& & \hspace{-0.5cm}  \text{s.t.}& \quad \omicron[\boldsymbol{\omega}](\mathbf{d}^t) = (\mathbf{\dot{g}}^t) & \forall t \in \mathcal T 
\end{flalign}
and the Lagrangian dual is simply
\begin{flalign}
LD = \max_{\boldsymbol{\lambda}} \ LR(\boldsymbol{\lambda}) \label{abst.LD}
\end{flalign}  

\noindent
Problem \eqref{abst.LD} is solved by iterating between training for
weights $\boldsymbol{\omega}$ and updating the Lagrangian
multipliers. Iteration $j$ uses Lagrangian multiplier
$\boldsymbol{\lambda}^j$ and solves $LR(\boldsymbol{\lambda}^j)$ to
obtain the optimal weights $\boldsymbol{\omega}^j$.  It then updates
the Lagrangian multipliers using the constraint violations.  The
overall scheme is presented in Algorithm \ref{Alg.learn}. Lines
\ref{start-LR}--\ref{end-LR} train weights $\boldsymbol{\omega}^j$ for
a fixed vector of Lagrangian multipliers $\boldsymbol{\lambda}^j$,
using minibatches and a stochastic gradient descent method with
learning rate $\alpha$. For each minibatch, the algorithm computes the
predictions (line \ref{LR-Pred}), the constraint violations (line
\ref{LR-violations}), and updates the weights (line
\ref{LR-update}). Lines \ref{start-LD}--\ref{end-LD} describe the
solving of Lagrangian dual. It computes the Lagrangian relaxation
described previously and updates the Lagrangian multipliers in line
\ref{LD-update} using the median violation $\tilde{\nu}_c$ for each
nominal constraint $c$.

\begin{algorithm}[t]
	{\footnotesize \caption{{\small Lagrangian Dual Model ($\mathcal{T}$, $\mathbb{C}$, $\alpha$, $\rho$, {\it Jmax}, $\boldsymbol{\lambda}^0$, $\boldsymbol{\omega}^0$)}}\label{Alg.learn}
		\begin{algorithmic}[1]
			\State{$j\leftarrow0$.} 
			\For{$j=0, 1, \ldots, {\it Jmax}$} \label{start-LD}
			    \For{$k=0, 1, \ldots$} \label{start-LR}
			    \State{Sample minibatch: $\mathbb{T}_k\subset\mathcal{T}$}
				\For{$t\in\mathbb{T}_k$}
					\State{Compute \omicron$[\boldsymbol{\omega}^j](\mathbf{d}^t)= \mathbf{\dot{g}}^t $ and $\mathbb{L}^t_0(\mathbf{\dot{g}}^t)$  } \label{LR-Pred}
					\State{Compute $\nu^t_{c}(\mathbf{\dot{g}}^t)$ $\forall c\in\mathbb{C}$  }  \label{LR-violations}
				\EndFor
				\State{$\boldsymbol{\omega}^{j}\leftarrow \boldsymbol{\omega}^{j} - \alpha \nabla_{\boldsymbol{\omega}^{j}} [\sum_{t\in\mathcal{T}} (\mathbb{L}^t_0(\mathbf{\dot{g}}^t)  + \sum_{c\in\mathcal{C}} \lambda_c \nu^t_{c}(\mathbf{\dot{g}}^t))]$} \label{LR-update}
				\EndFor \label{end-LR}
				\State{$\lambda^{j+1}_c\leftarrow \lambda^{j}_c + \rho \, \tilde{\nu}_c \; \forall c\in \boldsymbol{C}$} \label{LD-update}
				\State{$\boldsymbol{\omega}^{j+1}\leftarrow \boldsymbol{\omega}^{j}$}
			\EndFor \label{end-LD}
		\end{algorithmic}
	} 
\end{algorithm}

\subsection{CCGA-DNN Model}
\label{sec.combining}

This section presents the final ML model, the CCGA-DNN, which mimics a
CCGA algorithm. In particular, the CCGA-DNN combines the Lagrangian
dual model with an outer loop that adds constraints for the contingent
states on-demand.

Observe first that a direct Lagrangian dual approach to the SCOPF
would require an outer loop to add predictors
$[\mathbf{\dot{g}}_{s},\mathbf{\dot{x}}_s,
  \dot{n}_s]_{s\in\mathcal{S}}$ and constraints
\eqref{ML.const.set02}--\eqref{ML.const.set03} for selected
contingency states $s$. Unfortunately, the addition of new predictors
structurally modifies the DNN output
\omicron$[\boldsymbol{\omega}](\cdot)$ and induces a considerable
increase in the DNN size.

The key idea to overcome this difficulty is to mimic the CCGA closely,
replacing the master problem with the prediction
$\omicron[\boldsymbol{\omega}^l](\mathbf{d}^t)$ at iteration $l$.
Moreover, constraints \eqref{ML.const.set03} are replaced by
constraints of the form 
\begin{align}
\dot{g}^t_{s,i} = \max\{0, \min\{\dot{g}^t_i + \dot{n}^t_s \gamma_i \,\hat{g}_i\,,\overline{g}_i\}\}, \label{cnt.approx.gs}
\end{align}
where $\dot{n}^t_s$ is obtained by the bisection method on the
prediction. Again, these constraints are not differentiable but admit
subgradients and hence can be dualized in the objective function.  The
CCGA-DNN is summarized in Algorithm \ref{Alg.CCGA_DNN}. At each
iteration $l$, the Lagrangian dual model (Algorithm \ref{Alg.learn})
produces updated weights $\boldsymbol{\omega}^l$ and multipliers
$\boldsymbol{\lambda}^l$ (line 5). The inner loop (lines 6--12)
applies the bisection method to find $\dot{n}^t_s$ for all $t$ and
constraints \eqref{disj01}--\eqref{disj06} to obtain
$\mathbf{\dot{g}}^t_s$ (lines 8--9). These values are then used to
compute the highest line violation $\phi^t$ among all states and the
associated contingent state $s_{\phi}^t$ (line 10). The inner loop
also increases the element of the counter vector $\mathbf{p}$
associated with $s_{\phi}^t$ whenever the highest violation for solve
$t$ is above tolerance $\epsilon$ (line 11). Then, in the main loop,
contingency states with high frequencies of violated lines are
identified (line 13) using a threshold $\beta_1$. The algorithm is
terminated if $\mathbb{S}^{'}$ is empty and median relative violations
for nominal constraints in \eqref{ML.const.set01} are within
tolerances $\beta_c$ (line 14). Otherwise, the set of constraints is
updated by adding constraints \eqref{Cut_01}--\eqref{Cut_02} and
\eqref{cnt.approx.gs} for added contingent states (line 15). The
Lagrangian multipliers for added constraints
\eqref{Cut_01}--\eqref{Cut_02} $\forall s\in \mathbb{S}^{'}$ are
initialized in line 16. Finally, constraints
\eqref{Cut_01}--\eqref{Cut_02} $\forall s\in \mathbb{S}$ are updated
with the median violation $\tilde{\phi}$ for those lines associated
with some $\phi^t$ (line 17). Note that the process of updating
Lagrangian multipliers for \eqref{Cut_01}--\eqref{Cut_02} $\forall
s\in \mathbb{S}$ is different and much stricter than that for nominal
constraints in Algorithm \ref{Alg.learn}.
\begin{algorithm}[t]
	{\footnotesize \caption{{\small CCGA-DNN ($\mathcal{T}$, $\alpha$, $\rho$, $\beta_1$, $\beta_c$, $\epsilon$, jMax) }}\label{Alg.CCGA_DNN}
		\begin{algorithmic}[1]
			\State{$\mathbb{C}\leftarrow\{\eqref{ML.const.set01}\}$, $\mathbb{S}\leftarrow\emptyset$,
			 $\boldsymbol{\lambda}^0\leftarrow \boldsymbol{0}$, $\boldsymbol{\omega}^0\leftarrow \boldsymbol{0}$}
			\State{Create a counter vector $\mathbf{p}$ of size $|\mathcal{S}|$ } 
			\For{$l=1,2, \ldots$}
			\State{$\mathbf{p}\leftarrow\mathbf{0}$}
			\State{$\boldsymbol{\lambda}^l$,\,$\boldsymbol{\omega}^l\leftarrow$ Lagrangian Dual Model($\mathcal{T}$, $\mathbb{C}$, $\alpha$, $\rho$, jMax, $\boldsymbol{\lambda}^{l-1}$,$\boldsymbol{\omega}^{l-1}$)}
			\For{$t\in\mathcal{T}$}
			    \State{$\mathbf{\dot{g}}^{t}\leftarrow$ \omicron$[\boldsymbol{\omega}^l](\mathbf{d}^t)$}
			    \State{$\dot{n}^t_s \leftarrow$ bisection method, $\forall s\in\mathcal{S}$}
			    \State{$\mathbf{\dot{g}}^{t}_s\leftarrow$ enforce \eqref{disj01}--\eqref{disj06}, $\forall s\in\mathcal{S}$}
			    \State{Compute: $\phi^t$ and identify  $s_{\phi}^t$  }
			    \State{\textbf{if} $\phi^t>\epsilon$ \textbf{then} increase $(s_{\phi}^t)$-th element of $\mathbf{p}$ by 1} 
				\EndFor
			\State{$\mathbb{S}^{'} \leftarrow \{s \,| \,\mathbf{p}[s]\,/\,|\mathcal{T}| > \beta_1$\}}
			\State{\textbf{BREAK} if $\mathbb{S}^{'}\equiv \emptyset$ and $\tilde{\nu_c}\leq\beta_c, \forall c\in \eqref{ML.const.set01}$.}
			\State{$\mathbb{C}\leftarrow\mathbb{C}\cup \{\eqref{Cut_01}-\eqref{Cut_02}, \eqref{cnt.approx.gs}, \, \forall s\in(\mathbb{S}^{'}\setminus\mathbb{S})\} $}
			\State{$\lambda^{l}_{\eqref{Cut_01}},\,\lambda^{l}_{\eqref{Cut_02}}\leftarrow 0,\; \forall s\in(\mathbb{S}^{'}\setminus\mathbb{S})$}
			\State{$\lambda^{l}_{\eqref{Cut_01}},\,\lambda^{l}_{\eqref{Cut_02}}\mathrel{+}= \rho \, \tilde{\phi},\; \forall s\in\mathbb{S}$} 
			\EndFor
		\end{algorithmic}
	} 
	\vspace{-0.1cm}
\end{algorithm}      


\section{Feasibility Recovery and Optimality Gap}
\label{sec.restaur}  
\label{sec.optimality}

The training step produces a set of weights $\boldsymbol{\omega}^*$
and the associated DNN produces, almost instantly, a dispatch
prediction $\mathbf{\dot{g}} =
\omicron[\boldsymbol{\omega}^*](\mathbf{d})$ for an input load vector
$\mathbf{d}$. However, the prediction $\mathbf{\dot{g}}$ may violate
the nominal and contingency constraints. To restore feasibility, this
paper proposes a feasibility-recovery CCGA, denoted by FR-CCGA, that
finds the feasible solution closest to $\mathbf{\dot{g}}$. The master
problem for FR-CCGA is similar to
\eqref{master.ccga}--\eqref{Cut_02.ccga} but it uses a different
objective function, i.e.,
\begin{align}
\min_{\mathbf{g},[\mathbf{g}_{s}^{'}]_{s\in\mathcal{S}},[\mathbf{x}_{s},n_s]_{s\in\mathbb{S}}}
\quad&   |\hspace{-0.01cm}|\mathbf{\dot{g}} - \mathbf{g}|\hspace{-0.03cm}|\label{obj.FRCCGA} \\
\text{s.t.:}\quad & \eqref{const.nom}\text{--}\eqref{Cut_02.ccga} & \forall \; \mathcal{S},\, \mathbb{S},\, \mathbb{U}^+, \,\mathbb{U}^-\label{const.FRCCGA}
\end{align}    
Note that $\mathbf{\dot{g}}$ is a constant vector in FR-CCGA. While
CCGA and FR-CCGA are similar in nature, FR-CCGA is significantly
faster because $\omicron[\boldsymbol{\omega}^*](\mathbf{d})$ is often
close to feasibility.

The FR-CCGA and CCGA can be run in parallel to provide upper and lower
bounds to the SCOPF respectively. This may be useful for operators to
assess the quality of the prediction and the associated FR-CCGA solution
and decide whether to commit to the FR-CCGA solutions or wait until a better
solution is found or the optimality gap is sufficiently small. 
  
\section{Computational Experiments}
\label{sec.results}

\subsection{Data}
\label{sec.data}

The test cases are based on modified versions of 3 system topologies
from \cite{PGLIB}. Table \ref{tab:Inst.size} shows the size of a
single instance for each topology. For each topology, the training and
testing data is given by the inputs and solutions of many instances
that are constructed as follows. For each instance, the net demand of
each bus has a deterministic component and a random component. The
deterministic component varies across instances from 82\% of the
nominal net load to near-infeasibility values by small increments of
0.002\%. The random component is independently and uniformly
distributed ranging from -0.5\% to 0.5\% of the corresponding nominal
nodal net load for each bus and instance. Algorithm
\ref{Alg.CCGA_feas} was applied to solve each instance for a maximum
line violation of $\epsilon$ = 0.05 MW and an optimality gap of
0.25\%.

\begin{table}[!t]
	\centering\scriptsize
	\caption{Instance Size for the SCOPF Problem \eqref{obj.func}--\eqref{const.set03} after  Presolve}
	\vspace{-0.28cm}
	\begin{tabular}{l r r r r r r}
		\toprule
		\noalign{\vskip 0.35mm}
	    \multirow{2}{0.35cm}{System} & \multirow{2}{0.32cm}{$|\mathcal{G}|\;\,$} & \multirow{2}{0.32cm}{$|\mathcal{L}|\;\,$} & \multirow{2}{0.32cm}{$|\mathcal{N}|\;\,$}	& \multicolumn{1}{c}{Total} &  \multicolumn{1}{c}{Binary} &  \multicolumn{1}{c}{Linear}  \\
								 & & &		&  \multicolumn{1}{c}{Variables} & \multicolumn{1}{c}{Variables}	&  \multicolumn{1}{c}{Constraints}    \\
		\midrule
		118-IEEE	   &54&186     &118&13,466&2,862&19,137\\
		1354-PEG    &260&1,991&1,354&387,026&63,455&513,677\\
		1888-RTE       &297&2,531&1,888&467,011&79,032&624,780\\
		\bottomrule
			\vspace{-0.35cm}
	\end{tabular}%
	\label{tab:Inst.size}%
\end{table}
\vspace{-0.1cm}
\subsection{Training Aspects}

The training set $\mathcal{T}$ is composed by a random sample
containing 70\% of the generated instances. Algorithm
\ref{Alg.CCGA_DNN} was applied with $\epsilon$ set to 1 MW, $\beta_1$
to 5\%, $\beta_c$ to $1.5\cdot 10^{-2}$, and $\rho$ to $10^5$.  The
inner loop of Algorithm \ref{Alg.learn} (lines 4--11) is executed $1.5
\cdot 10^5$ times with a learning rate $\alpha$ varying from $10^{-4}$
to $10^{-10}$ and {\em Jmax} was set to 1.  The DNN models were
implemented using PyTorch package with Python 3.0. The training was
performed using NVidia Tesla V100 GPUs and 2 GHz Intel Cores. Table
\ref{tab:Trainingsumm} presents a training summary.  In the following,
the baseline model is denoted by $\mathcal{M}_{b}$ and the
CCGA-DNN by $\mathcal{M}_{ccga}$.  Note that the first
iteration of Algorithm \ref{Alg.CCGA_DNN} returns the weights of the
baseline model.

\begin{table}[t!]
	\centering\scriptsize
	\caption{Training Summary for Algorithm \ref{Alg.CCGA_DNN}}
	\vspace{-0.28cm}
	\begin{tabular}{c c c }
		\toprule
		System & Iterations & Added Contingency States $\mathbb{S}$ (by generator numbers) 
		\\[0.3em]
		\midrule
		118-IEEE  &3                & \{4\}\\
		1354-PEG  &3                & \{23, 65, 74, 112, 126, 163, 222\} \\
		1888-RTE  &2                & \{152, 153\} \\
		\bottomrule
	\vspace{-0.4cm}
	\end{tabular}%
	\label{tab:Trainingsumm}%
\end{table}
\vspace{-0.1cm}
\subsection{Prediction Quality}

\begin{figure}[!t]
	\centering
	\includegraphics{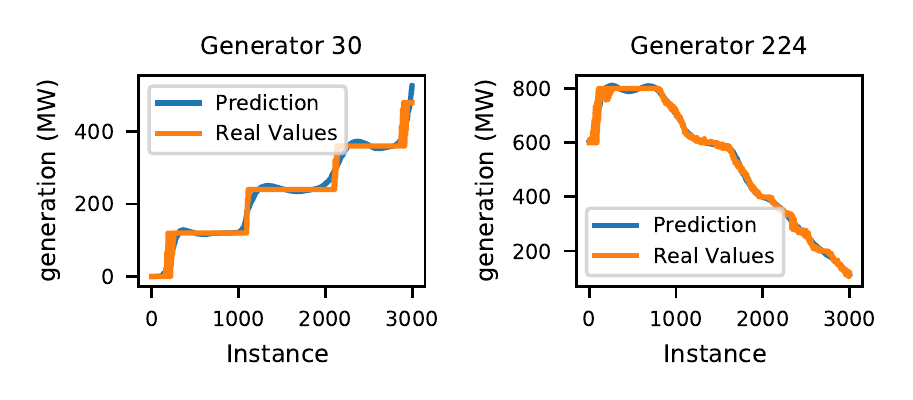}
	\vspace{-0.69cm}
	\caption{Prediction of $\mathcal{M}_{ccga}$ for selected generators of the 1354-PEG System.}
	\label{fig_sim}
\end{figure}

Accurate predictions were obtained for all DNN models
and topologies. Figure \ref{fig_sim} illustrates how
$\mathcal{M}_{ccga}$ can learn complex generator patterns arising in the
1354-PEG system. Table \ref{tab:Predict.err} reports the mean absolute
errors for predictions $\dot{\mathbf{g}}^t$, segmented by generation
range: $\mathcal{M}_{b}$ achieves a slightly better accuracy
which is expected since it is the less constrained
model. 
\begin{table}[t!]
	\centering\scriptsize
	\caption{Prediction Mean Absolute Errors (\%)}
	\vspace{-0.28cm}
	\begin{tabular}{c l r r r r r r r}
		\toprule
		& &\multicolumn{7}{c}{Generation Range (MW)}\\[0.1em]
		\cline{3-9}\\[-0.45em]
	System	&Model	&10& 50 &100&250&500&1000 &2000 \\
		&	&50&100&250&500&1000&2000 &5000 \\
		\midrule 
		\multirow{2}{*}{118-IEEE}&	$\mathcal{M}_{b}$             &2.3&2.8&0.7&0.3&N/A&N/A&N/A\\
		&$\mathcal{M}_{ccga}$             &2.5&3.0&0.7&0.4&N/A&N/A&N/A\\
		\midrule
		\multirow{2}{*}{1354-PEG}&$\mathcal{M}_{b}$              &2.4&1.3&1.1&0.9&0.4&0.2&0.1\\
		&$\mathcal{M}_{ccga}$             &5.0&1.8&1.2&1.0&0.4&0.3&0.2\\
		\midrule
	    \multirow{2}{*}{1888-RTE}&$\mathcal{M}_{b}$             &1.3&1.2&0.7&0.4&0.3&0.1&N/A\\
		&$\mathcal{M}_{ccga}$           &1.4&1.1&0.6&0.4&0.3&0.1&N/A\\
		\bottomrule
		\vspace{-0.35cm}
	\end{tabular}%
	\label{tab:Predict.err}%
\end{table}


Table \ref{tab:ConstraintViol} reports selected indicators of
violations: the relative violation
$\lambda_{\eqref{PTDF_nominal_first}}$ of the total load constraint and
the relative violation RLV of the lines associated with $\phi$. The
results report median values as well as lower and upper bounds for
intervals that capture 95\% of the instances. Both models achieve the
desired tolerance of $\beta_c = 1.5\cdot10^{-2}$ for
$\lambda_\eqref{PTDF_nominal_first}$ (the tolerance $\beta_c$ does not
apply to RLV). Model $\mathcal{M}_{ccga}$ produces lower overall
violations and has a major effect on RLV.

\begin{table}[t!]
	\centering\scriptsize
	\caption{Selected Indicators of Violation Across Instances (\%)}
	\vspace{-0.28cm}
	\begin{tabular}{c l r r r | r r r }
		\toprule
		 & 	&\multicolumn{3}{c|}{$\lambda_\eqref{PTDF_nominal_first}$}&\multicolumn{3}{c}{RLV}\\[0.05cm]
		 \cline{3-5}\cline{6-8}\\[-0.1cm]
		System & Model &\multicolumn{1}{c}{Median}&\multicolumn{2}{c|}{95\%-Interval}&\multicolumn{1}{c}{Median}&\multicolumn{2}{c}{95\%-Interval} \\
		\midrule
		\multirow{2}{*}{118-IEEE}& $\mathcal{M}_{b}$	            &0.016&0.001&0.055&0.099&0.000&0.467\\ 
		                            &$\mathcal{M}_{ccga}$              &0.003&0.000&0.011&0.000&0.000&0.548\\ 
		\midrule
		\multirow{2}{*}{1354-PEG}&$\mathcal{M}_{b}$       &0.018&0.001&0.062&0.259&0.086&1.256\\ 
		&$\mathcal{M}_{ccga}$   &0.010&0.000&0.040&0.005&0.000&0.117\\ 
		\midrule
		\multirow{2}{*}{1888-RTE}&$\mathcal{M}_{b}$      &0.012&0.001&0.047&0.205&0.021&1.418\\ 
	    &$\mathcal{M}_{ccga}$ &0.005&0.000&0.027&0.007&0.000&0.057\\ 
		\bottomrule
		\noalign{\vskip 0.38mm}
		\multicolumn{8}{l}{$\lambda_\eqref{PTDF_nominal_first}$ -- Net load constraint violation divided by total load.}\\
		\multicolumn{8}{l}{RLV -- Relative violation for line associated with $\phi$.}\\
	\vspace{-0.35cm}
	\end{tabular}%
	\label{tab:ConstraintViol}%
\end{table}

\vspace{-0.1cm}
\subsection{Comparison with Benchmark CCGA}

The previous sections reported on the accuracy of the predictors. This
section shows how FR-CCGA leverages the predictors to find
near-optimal primal solutions significantly faster than CCGA.  More
precisely, it compares, in terms of cost and CPU time, CCGA and
FR-CCGA when seeded with $\mathcal{M}_{b}$ and $\mathcal{M}_{ccga}$,
for $200$ randomly selected instances for each system topology. Each
instance was solved with the same tolerances as in the training. They
were solved using Gurobi 8.1.1 under JuMP package for Julia 0.6.4 on a
laptop Dell XPS 13 9380 featuring a i7-8565U processor at 1.8 GHz and
16 GB of RAM. Tables \ref{tab:distance}, \ref{tab:CorrecFeasibility},
and \ref{tab:Correct_Cost_Increase} summarize the experiments.  Table
\ref{tab:distance} reports the distances in percentage between the
predictions and the feasible solutions obtained by FR-CCGA when seeded
with the predictions. For instance, for $\mathcal{M}_{b}$, this
distance is $\sum_i | g_i^b - g_i^f | / \sum_i g_i^f$, where
$g_i^b$ is $\mathcal{M}_{b}$'s prediction for generator $i$ and
$g_i^f$ is  $i$'s generation computed by FR-CCGA seeded with
$\boldsymbol{g}^b$.  As should be clear, the predictions are very
close to feasibility. Table \ref{tab:CorrecFeasibility} reports the
computation times which show significant increases in performance by
FR-CCGA especially when seeded with $\mathcal{M}_{ccga}$ and on the
1354-PEG system, the most challenging network. FR-CCGA is about 160
times faster than CCGA on this test case. FR-CCGA is also
significantly more robust when using $\mathcal{M}_{ccga}$ instead of
$\mathcal{M}_{b}$. Table \ref{tab:Correct_Cost_Increase} indicates
that the cost/objective increase of FR-CCGA over CCGA is very small for both
$\mathcal{M}_{b}$ and $\mathcal{M}_{ccga}$.

\begin{table}[t!]
	\centering\scriptsize
	\caption{Distance between Prediction and Feasible Solution (\%)}
	\vspace{-0.28cm}
	\begin{tabular}{c l  r r r}
		\toprule
		System & Model & \multicolumn{1}{c}{Median} & \multicolumn{2}{c}{95\%-Interval} 	 \\[0.3em]
		\midrule
		\multirow{2}{*}{118-IEEE}&$\mathcal{M}_{b}$	            & 0.05&0.01&0.12\\
		&$\mathcal{M}_{ccga}$	        &0.18&0.05&0.47\\
		\midrule
		\multirow{2}{*}{1354-PEG}&$\mathcal{M}_{b}$           &0.05&0.01&0.13\\
		&$\mathcal{M}_{ccga}$        &0.09&0.05&0.14\\
		\midrule
		\multirow{2}{*}{1888-RTE}&$\mathcal{M}_{b}$                &0.04&0.00&0.10\\
		&$\mathcal{M}_{ccga}$            &0.04&0.01&0.07\\
		\bottomrule
		\vspace{-0.4cm}
	\end{tabular}%
	\label{tab:distance}%
\end{table}


\begin{table}[t!]
	\centering\scriptsize
	\caption{CPU Time Comparison}
	\vspace{-0.28cm}
	\begin{tabular}{c l r r r r r }
		\toprule
		System & Model &\multicolumn{1}{c}{Median}&\multicolumn{1}{c}{Mean}&\multicolumn{1}{c}{Min.}&\multicolumn{1}{c}{Max.}&\multicolumn{1}{c}{Std.} \\
		\midrule
		\multirow{3}{*}{118-IEEE}&CCGA & 0.210&0.214&0.101&3.717&0.305\\
		            & $\mathcal{M}_{b}$ &0.024&0.057&0.021&1.580&0.160\\
	                &$\mathcal{M}_{ccga}$   &0.026&0.068&0.023&1.372&0.171\\ 
    		\midrule
		\multirow{3}{*}{1354-PEG}&CCGA  &321.746&327.210&75.585&741.798&127.101\\	                
	                            &$\mathcal{M}_{b}$   &5.335&8.434&1.320&133.366&13.505\\
	                            &$\mathcal{M}_{ccga}$  &1.521&2.168&0.768&8.449&1.740\\  
	      \midrule
	   	\multirow{3}{*}{1888-RTE}&CCGA  &5.479&7.406&3.110&30.923&7.073\\
	                         &$\mathcal{M}_{b}$ &5.316&5.501&1.224&18.074&3.348\\ 
	                         &$\mathcal{M}_{ccga}$  &2.120&1.945&0.911&4.543&0.919\\ 
		\bottomrule
		\noalign{\vskip 0.38mm}
	\end{tabular}%
	\vspace{-0.35cm}
	\label{tab:CorrecFeasibility}%
\end{table}

\begin{table}[t!]
	\centering\scriptsize
	\caption{FR-CCGA Cost Increase over CCGA (\%)}
	\vspace{-0.30cm}
	\begin{tabular}{c l r r r r r }
		\toprule
		System & Model &\multicolumn{1}{c}{Median}&\multicolumn{1}{c}{Mean}&\multicolumn{1}{c}{Min.}&\multicolumn{1}{c}{Max.}&\multicolumn{1}{c}{Std.} \\
		\midrule
		\multirow{2}{*}{118-IEEE}&$\mathcal{M}_{b}$ &0.021&0.019&-0.073&0.055&0.014\\
		                         &$\mathcal{M}_{ccga}$&0.027&0.030&-0.010&0.112&0.020\\ 
    		\midrule
		\multirow{2}{*}{1354-PEG}&$\mathcal{M}_{b}$   &0.020&0.021&-0.007&0.051&0.012\\
	                            &$\mathcal{M}_{ccga}$   &0.067&0.067&0.032 &0.091&0.017\\  
	      \midrule
	   	\multirow{2}{*}{1888-RTE}&$\mathcal{M}_{b}$ &0.026&0.024&-0.003&0.070&0.011\\ 
	                         &$\mathcal{M}_{ccga}$  &0.033&0.033&0.004&0.070&0.013\\ 
		\bottomrule
		\noalign{\vskip 0.38mm}
	\end{tabular}%
	\label{tab:Correct_Cost_Increase}%
	\vspace{-0.35cm}
\end{table}

\begin{figure}[!t]
	\centering
	\includegraphics{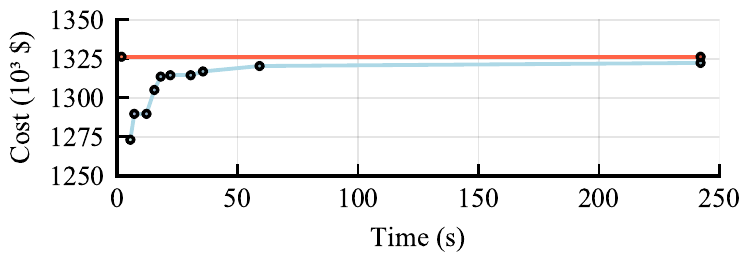}
	\vspace{-0.3cm}
	\caption{Convergence plot for the $\mathcal{M}_{ccga}$ for the 1354-PEG system.}
	\label{fig_sim2}
\end{figure}

Figure \ref{fig_sim2} illustrates the behavior of the algorithms on a
randomly chosen instance of the 1354-PEG system. The red line
represents the upper bound (feasible solution) generated in $1.87$
seconds by the FR-CCGA seeded with $\mathcal{M}_{ccga}$. The blue line
represents a sequence of true lower bounds (infeasible solutions)
generated by Algorithm \ref{Alg.CCGA_feas}.
\vspace{-0.1cm}

\section{Conclusion}
\label{sec.conclusion}

This paper proposed a tractable methodology that combines deep
learning models and robust optimization for generating solutions for
the SCOPF problem. The considered SCOPF modeled generator
contingencies and the automatic primary response of synchronized
units. Computational results over two large test cases demonstrate the
practical relevance of the methodology as a scalable, easy to specify,
and cost-efficient alternative tool for managing short-term
scheduling.

\end{document}